\documentstyle{amsart}
\def\dj{d\kern-.30em\raise1.25ex\vbox{\hrule width .3em height .03em}}
\def\Dj{D\kern-.70em\raise0.75ex\vbox{\hrule width .3em height .03em}
\kern.03em}
\newsymbol\restr 1316
\newcommand{\cl}{\mbox{\family{euf}\shape{n}\selectfont cl}}
\newcommand{\e}{\epsilon}
\newcommand{\k}{\kappa}
\newcommand{\id}{\mbox{\shape{n}\selectfont id}}
\newcommand{\Sum}{\displaystyle{\sum}}
\newcommand{\im}{\mbox{im}}
\newcommand{\ita}{\iota^F}
\newcommand{\imb}{\ell}
\newcommand{\lin}{\mbox{lin}}
\newtheorem{lem}{Lemma}

\begin{document}
\title[Braided Clifford Algebras]
{Braided Clifford Algebras as Braided\\ Quantum Groups}
\author{Mi\'co \Dj ur\Dj evi\'c}
\address{Instituto de Matematicas, UNAM, Area de la Investigacion
Cientifica, Circuito Exterior, Ciudad Universitaria, M\'exico DF,
CP 04510, MEXICO\newline
\indent{\it Written In}\newline
\indent Facultad de Estudios Superiores Cuautitlan, UNAM, MEXICO}
\maketitle
\begin{abstract}
The paper deals with braided Clifford algebras, understood as
Chevalley-K\"ahler deformations of braided exterior algebras. It is
shown that Clifford algebras based on involutive braids
can be naturally endowed
with a braided quantum group structure. Basic group entities are
constructed explicitly.
\end{abstract}
\section{Introduction}
The aim of this letter is to present a general construction of examples
of braided quantum groups, based on Clifford algebras associated
to involutive braids. From the
noncommutative-geometric \cite{C} point of view, these structures include
completely pointless objects, and overcome in such a way an inherent
geometrical inhomogeneity of standard quantum groups (caused
by the presence of classical points).
{\renewcommand\thepage{}
The construction is based on a general theory of
braided Clifford algebras \cite{MZ}, interpreted as Chevalley-K\"ahler
deformations of the corresponding braided exterior algebras. It will be
shown that if the braid operator is involutive the corresponding
Clifford algebras can be
naturally equipped with a coalgebra structure and the antipode map, so
that the whole structure becomes a braided quantum group \cite{D}.

The paper is organized as follows. In the next two sections basic
properties of
braided quantum groups and braided Clifford algebras are
collected. The main construction of the group structure on braided
Clifford algebras is given in
Section~4. The construction is presented in two conceptually different
ways. Finally, in the last section some concluding remarks are made.
\section{Braided Exterior and Clifford Algebras}
Let $W$ be a (complex) finite-dimensional
vector space, and let $\psi$
be an arbitrary linear automorphism of $W\otimes W$
satisfying the braid equation
$$
(\psi\otimes\id)(\id\otimes\psi)(\psi\otimes\id)
=(\id\otimes\psi)(\psi\otimes\id)(\id\otimes\psi).
$$

By definition \cite{W},
the corresponding {\it exterior algebra} $W^\wedge$ is the
factoralgebra of the tensor algebra $W^\otimes$ relative
to the ideal $\ker(A)\subseteq W^\otimes$.
Here $A\colon W^\otimes\rightarrow W^\otimes$
is the corresponding total antisymetrizer map. Its components $A_n\colon
W^{\otimes n}\rightarrow W^{\otimes n}$ are given by
$$ A_n=\sum_{\pi\in S_n} (-1)^\pi\psi_\pi$$
where $\psi_\pi\colon W^{\otimes n}\rightarrow W^{\otimes n}$ are
maps obtained by replacing transpositions
figuring in a minimal decomposition of  $\pi$
by the corresponding $\psi$-twists. The following factorizations hold
\begin{equation} A_{n+k}=(A_n\otimes A_k)A_{nk}\qquad
A_{n+k}=B_{nk}(A_n\otimes A_k)\label{dec2}
\end{equation}
where
$$
A_{nk}=\sum_{\pi\in S_{nk}}(-1)^\pi\psi_{\pi^{-1}}\qquad
B_{nk}=\sum_{\pi\in S_{nk}}(-1)^\pi\psi_\pi
$$
and $S_{nk}\subseteq S_{n+k}$ is the set of permutations preserving the
order of sets $\{1,\dots,n\}$ and $\{n+1,\dots,n+k\}$.

The algebra $W^\wedge$ can be naturally realized as a subspace
$\im(A)\subseteq W^\otimes$. This realization is given by
$$ \Bigl[\zeta+\ker(A)\Bigr]\, \leftrightarrow\, A(\zeta).$$
In terms of the above identification,
$$\zeta\wedge\xi=B_{nk}(\zeta\otimes\xi)$$
for each $\zeta\in W^{\wedge n}$ and $\xi\in W^{\wedge k}$.

The map $\psi$ can be naturally extended to a braiding on
$W^{\wedge}$ (which will be denoted by the same symbol) such that
\begin{align}
\psi(\id\otimes m)&=(m\otimes\id)(\id\otimes\psi)(\psi\otimes
\id)\label{ps-m1}\\
\psi(m\otimes\id)&=(\id\otimes m)(\psi\otimes\id)(\id\otimes
\psi)\label{ps-m2},
\end{align}
where $m$ is the corresponding product map. By construction,
$$\psi(W^{\wedge i}\otimes W^{\wedge j})=W^{\wedge j}\otimes W^{\wedge
i}$$ for each $i,j\geq 0$.

Let us assume that $W$ is endowed with a ``quadratic form'' $F$, which will
be understood
as a linear map $F\colon W\otimes W\rightarrow \Bbb{C}$ (or equivalently,
as a bilinear map $F\colon W\times W\rightarrow \Bbb{C}$).
Further, let us
assume that $F$ and $\psi$ are mutually related in the following way
\begin{equation}\label{ffun}
(F\otimes \id)(\id\otimes\psi)=(\id\otimes F)(\psi\otimes\id)
\end{equation}
}
Let $\ita\colon W\times W^\wedge\rightarrow W^\wedge$ be a bilinear
map specified by
\begin{align*}
\ita_x(\zeta)&=F(x,\zeta)\\
\ita_x(\vartheta\wedge\eta)&=\ita_x(\vartheta)\wedge
\eta+(-1)^{\partial\vartheta}
\sum_k\vartheta_k\wedge\ita_{x_k}(\eta)
\end{align*}
where $\Sum_k\vartheta_k\otimes x_k=\psi(x\otimes\vartheta)$ and
$\zeta\in W$.
The second condition is a braided variant of the Leibniz rule.
The map $\ita$ is uniquely
(and consistently) determined by the above conditions. We have
$$\ita_x(1)=0$$
for each $x\in W$.

The introduced contraction operator can be trivially extended to the map
of the form $\ita\colon W^{\otimes}\times W^{\wedge}\rightarrow W^{\wedge}$
by requiring
$$\ita_{u\otimes w}=\ita_u\ita_v$$
for each $u,v\in W^\otimes$. It turns out that
if $u\in\ker(A)$ then $\ita_u=0$.

Therefore, it is possible to pass from
$W^\otimes$ to $W^\wedge$ in the first argument
of $\ita$. In such a way we obtain a contraction map $\ita
\colon W^\wedge\times W^\wedge\rightarrow W^\wedge$.

Let us define  ``relative'' contraction operators
$\langle\,\rangle_k\colon W^\wedge\times W^\wedge\rightarrow W^\wedge$ as
bilinear maps of the form
$$\langle\vartheta,\eta\rangle_k=
\sum_{\alpha\beta}\alpha\wedge\ita_\beta(\eta)$$
where $$\sum_{\alpha\beta}\alpha\wedge\beta=A_{lk}(\vartheta^*)$$
while $\alpha\in W^{\wedge l}$, $\beta\in W^{\wedge k}$ and
$\vartheta\in W^{\wedge n}$, with $n=k+l$. Finally, $\vartheta^*\in
W^{\otimes n}$ is such that $[\vartheta^*]^\wedge=\vartheta$.
Consistency of this definition follows from the first decomposition in
\eqref{dec2}.

The formula
\begin{equation}\label{prodcli}
\vartheta\circ\eta=\vartheta\wedge\eta+\sum_{k\geq 1}
\langle\vartheta,\eta\rangle_k
\end{equation}
defines a new $F$-dependent product in the space $W^\wedge$. The
associativity of this product follows from \eqref{ffun}.
In particular,
$$x\circ\vartheta=x\wedge\vartheta+\ita_x(\vartheta)$$
for $x\in W$. This is a counterpart of classical Chevalley's formula.

Endowed with $\circ$, the space $W^\wedge$ becomes
a unital associative algebra, with the unity $1\in W^\wedge$.

By definition, the algebra
$\cl(W,\psi,F)=(W^\wedge,\circ)$ is called {\it
the braided Clifford algebra} (associated to $\{\psi,F\}$).

The constructed algebra can be understood as a deformation of
the exterior algebra $W^\wedge$. Furthermore, the graded algebra associated
to the filtered $\cl(W,\psi,F)$ naturally coincides with the braided exterior
algebra $W^\wedge$.

The construction of braided Clifford algebras can be performed in a
conceptually different way \cite{MZ},
introducing  a Clifford product in the tensor algebra $W^\otimes$.
Let us consider a map $\lambda_F\colon W^\otimes\rightarrow
W^\otimes$ specified by
$$\lambda_F(1)=1\qquad\lambda_F(x\otimes\vartheta)=x\otimes
\lambda_F(\vartheta)+\ita_x\lambda_F(\vartheta)$$
where $x\in W$ and $\vartheta\in W^\otimes$. In the above formula
$\ita$ is considered as a contraction acting in $W^\otimes$.
The map $\lambda_F$ is bijective. Let
$\circ$ be a new product in $W^\otimes$, given by
$$\vartheta\circ\eta=\lambda_F\bigl(\lambda_F^{-1}(\vartheta)
\otimes\lambda_F^{-1}(\eta)\bigr).$$

By construction the space $\ker(A)$ is a left ideal in $W^\otimes$,
relative to this new product. Condition \eqref{ffun} implies that
$\ker(A)$ is also a right $\circ$-ideal.

Let $J_F\subseteq W^\otimes$ be the ideal corresponding to the constructed
braided Clifford algebra. In other words,
$$\cl(W,\psi,F)=W^\otimes/ J_F,$$
in a natural manner. We have then
$$\lambda_F^{-1}\bigl[\ker(A)\bigr]=J_F\qquad(W^\otimes,\circ)/\ker(A)=
\cl(W,\psi,F).$$

\section{Basic Properties of Braided Quantum Groups}

Let us consider a complex associative algebra $\cal{A}$, with
the unit element $1\in\cal{A}$ and
the product map
$m\colon\cal{A}\otimes \cal{A}\rightarrow\cal{A}$. Further, let  us
assume that $\cal{A}$ is
equipped with a  coassociative  coalgebra
structure, specified by the coproduct
$\phi\colon\cal{A}\rightarrow\cal{A}\otimes
\cal{A}$ and the counit $\e\colon\cal{A}\rightarrow\Bbb{C}$.
Finally, let us assume that there exist
bijective linear maps $\k\colon\cal{A}\rightarrow\cal{A}$ and
$\sigma\colon\cal{A}\otimes
\cal{A}\rightarrow\cal{A}\otimes \cal{A}$ such that  the  following
equalities hold:
\begin{gather}
\sigma(m\otimes\id)=(\id\otimes m)(\sigma\otimes\id)
(\id\otimes \sigma)\label{1}\\
\sigma(\id\otimes m)=(m\otimes\id)(\id\otimes\sigma)(\sigma\otimes
\id)\label{2}\\
1\e=m(\id\otimes\k)\phi=m(\k\otimes\id)\phi\label{5}\\
\phi m=(m\otimes m)(\id\otimes\sigma\otimes\id)(\phi\otimes\phi)\label{3}\\
(\sigma\otimes\id^2)(\id\otimes\phi\otimes\id)(\sigma^{-1}\otimes\id)
(\id\otimes\phi)\phantom{(\id\otimes\phi\otimes\id)(\sigma^{-1}\otimes\id)}
\label{4}\\
\phantom{(\id\otimes\phi\otimes\id)(\sigma^{-1}\otimes\id)}
=(\id^2\otimes\sigma)(\id\otimes\phi
\otimes\id)(\id\otimes
\sigma^{-1})(\phi\otimes\id).\notag
\end{gather}

By definition \cite{D}, every pair
$G=\bigl(\cal{A},\bigl\{\phi,\e,\k,\sigma\bigr\}\bigr)$
satisfying the above requirements is called
{\it a braided quantum group}.

The map $\sigma$ plays the role of the ``twisting  operator'' (the
ordinary transposition in
the standard theory of Hopf algebras).
This operator induces a product in $\cal{A}\otimes\cal{A}$, via the
formula
$$
 (a\otimes b)(c\otimes d)=a\sigma(b\otimes c)d.
$$
Identities \eqref{1}--\eqref{2} ensure that this defines
an  associative  algebra
structure on $\cal{A}\otimes \cal{A}$, such that $1\otimes 1$  is  the
unit element. In particular,
$$
\sigma(1\otimes a)=a\otimes 1 \qquad\quad
\sigma(a\otimes 1)=1\otimes a.
$$
for each $a\in\cal{A}$. Equality \eqref{3}  then
says that $\phi$ is multiplicative. Identity \eqref{4}
expresses  the  coassociativity  of  the  map
$(\id\otimes \sigma^{-1}  \otimes \id)(\phi\otimes  \phi)$.

The
antipode $\k$ is uniquely determined by
\eqref{5}. The flip-over operator $\sigma$ is expressible through
$\phi,m$ and $\k$ in the following way
$$
\sigma=(m\otimes m)(\k\otimes\phi m\otimes\k)(\phi\otimes \phi).
$$

Besides the flip-over operator $\sigma$, a ``secondary''
flip-over operator $\tau$ naturally enters the game. This map is
specified by
$$
\tau\sigma^{-1}=(\id^2 \otimes \e)(\id\otimes\sigma^{-1})(\phi\otimes
\id)=(\e\otimes \id^2)(
\sigma^{-1}\otimes \id)(\id\otimes\phi).
$$
Its inverse is given by
$$
 \tau^{-1}\sigma=(\id^2 \otimes \e)(\id\otimes \sigma)(\phi\otimes
\id)=(\e\otimes\id^2)
(\sigma\otimes\id)(\id\otimes \phi).
$$

Two maps satisfy the following system of braid equations
\begin{align*}
(\sigma\otimes \id)(\id\otimes \sigma)(\sigma\otimes \id)&=(\id\otimes
\sigma)(\sigma\otimes \id)(\id\otimes \sigma)\\
(\tau\otimes \id)(\id\otimes \sigma)(\sigma\otimes \id)&=(\id\otimes \sigma)(
\sigma\otimes \id)(\id\otimes \tau)\\
(\sigma\otimes \id)(\id\otimes \tau)(\sigma\otimes \id)&=(\id\otimes \sigma
)(\tau\otimes \id)(\id\otimes \sigma)\\
(\sigma\otimes \id)(\id\otimes \sigma)(\tau\otimes \id)&=(\id\otimes \tau)(
\sigma\otimes \id)(\id\otimes \sigma)\\
(\tau\otimes \id)(\id\otimes \tau)(\sigma\otimes\id)
&=(\id\otimes \sigma)(\tau\otimes \id)(\id\otimes \tau)\\
(\tau\otimes \id)(\id\otimes \sigma)(\tau\otimes
\id)&=(\id\otimes \tau)(\sigma\otimes \id)(\id\otimes \tau)\\
(\sigma\otimes \id)(\id\otimes \tau)(\tau\otimes \id)&=(\id\otimes \tau)(
\tau\otimes \id)(\id\otimes \sigma)\\
(\tau\otimes \id)(\id\otimes \tau)(\tau\otimes \id)&=(\id\otimes
\tau)(\tau\otimes \id)(\id\otimes\tau).
\end{align*}

Operators $\sigma$ and $\tau$ naturally
generate a ``braid
system'' $\cal{T}=\bigl\{\sigma_n\mid n\in\Bbb{Z}\bigr\}$ consisting of
braidings of the form
$$\sigma_n=\tau(\sigma^{-1}\tau)^{-n}=(\tau\sigma^{-1})^{-n}\tau.$$

The following twisting properties hold
\begin{gather*}
(\phi\otimes\id)\sigma_{n+k}=(\id\otimes\sigma_k)(\sigma_n\otimes\id)
(\id\otimes\phi)\\
(\id\otimes\phi)\sigma_{n+k}=(\sigma_k\otimes\id)(\id\otimes\sigma_n)
(\phi\otimes\id)\\
\sigma_n(\k\otimes \id)=(\id\otimes \k)\sigma_{-n}\\
\sigma_n(\id\otimes \k)=(\k\otimes \id)\sigma_{-n}\\
\sigma_n(m\otimes\id)=(\id\otimes m)(\sigma_n\otimes\id)(\id\otimes\sigma_n)
\phantom{.}\\
\sigma_n(\id\otimes
m)=(m\otimes \id)(\id\otimes\sigma_n)(\sigma_n\otimes\id).
\end{gather*}

The formula
$$
\e m=(\e\otimes \e)\sigma^{-1}\tau
$$
generalizes the standard multiplicativity property for the counit.
 A  braided
counterpart of the standard anti(co)multiplicativity
property of the antipode is given by
\begin{align*}
\phi\k&=\sigma(\k\otimes \k)\phi\\
\k m&=m(\k\otimes \k)\tau\sigma^{-1}\tau\sigma^{-1}\tau.
\end{align*}

If $\sigma=\sigma^{-1}$ then all maps $\sigma_n$ are involutive, too.
\section{Construction of Group Structures}
This section is devoted to the construction of the canonical braided quantum
group structure on braided Clifford algebras. We shall assume that the
braid operator $\psi$ is {\it involutive}. Let $\cal{A}_F=\cl(W,\psi,F)$ be
the corresponding Clifford algebra.

The ideal $\ker(A)$ is generated by the space
$\ker(\id^2-\psi)=\im(\id^2+\psi)$ consisting of $\psi$-invariant
elements. It follows that
the corresponding Clifford ideal $J_F$ is generated by
elements of the form
\begin{equation}\label{gen}
Q=\sum_k\Bigl\{x_k\otimes y_k-F(x_k,y_k)1\Bigr\}
\end{equation}
where $\Sum_k\psi(x_k\otimes y_k)=\sum_kx_k\otimes y_k$.

We shall denote
by $m_F$ the Clifford product in $\cal{A}_F$. Without a lack of generality
we can assume that $F$ is $\psi$-symmetric.
\begin{lem}
There exists the unique linear
operator $\sigma_F\colon\cal{A}_F\otimes
\cal{A}_F\rightarrow\cal{A}_F\otimes\cal{A}_F$ satisfying
\begin{gather}
(m_F\otimes\id)(\id\otimes\sigma_F)(\sigma_F\otimes\id)=\sigma_F(\id\otimes
m_F)\label{m-S1}\\
(\id\otimes
m_F)(\sigma_F\otimes\id)(\id\otimes\sigma_F)=\sigma_F(m_F\otimes\id)
\label{m-S2}\\
\sigma_F(x\otimes y)=-\psi(x\otimes y)-F(x,y)1\otimes 1, \label{Sxy}
\end{gather}
for each $x,y\in W$. Furthermore, $\sigma_F$ is involutive and satisfies
the braid equation.
\end{lem}
\begin{pf}
Equalities of the same form
uniquely and consistently determine a linear operator
$\sigma_F\colon W^\otimes\otimes W^\otimes\rightarrow W^\otimes\otimes
W^\otimes$. The fact that $\sigma_F$ is a braid operator follows from
the braid equation for $\psi$, and the covariance property \eqref{ffun}.
Involutivity of $\psi$ and $\psi$-symmetricity of $F$ imply that
$\sigma_F$ is involutive.
Using the definition of $J_F$, and \eqref{ffun},
as well as the braid equation for $\psi$ it follows that
\begin{gather*}
\sigma_F(J_F\otimes W^\otimes)=W^\otimes\otimes J_F\\
\sigma_F(W^\otimes\otimes J_F)=J_F\otimes W^\otimes.
\end{gather*}
Hence, $\sigma_F$ can be projected on $\cal{A}_F\otimes\cal{A}_F$.
\end{pf}

In what follows it will be assumed that $\cal{A}_F\otimes\cal{A}_F$ is
endowed with the $\sigma_F$-induced product.
\begin{lem}
There exists the unique (unital) homomorphism $\phi\colon\cal{A}_F
\rightarrow\cal{A}_F\otimes\cal{A}_F$ satisfying
\begin{equation}\label{prim}
\phi(x)=1\otimes x+x\otimes 1
\end{equation}
for each $x\in W$.
\end{lem}
\begin{pf}
There exists the unique unital homomorphism $\phi\colon W^\otimes
\rightarrow\cal{A}_F\otimes\cal{A}_F$ satisfying the above requirement.
We have
$$\phi(J_F)=\{0\}.$$
Indeed it is sufficient to check that $\phi$
vanishes on elements of the form
\eqref{gen}. A direct calculation gives
\begin{multline*}
\phi(Q)=\sum_k\Bigl\{x_ky_k\otimes 1+x_k\otimes y_k+\sigma_F(x_k\otimes y_k)
+1\otimes x_ky_k-F(x_k,y_k)1\otimes 1\Bigr\}\\
=\sum_k\Bigl\{\bigl(x_ky_k-F(x_k,y_k)1\bigr)\otimes 1+1\otimes
\bigl(x_ky_k-F(x_k,y_k)1\bigr)\Bigr\}=0.
\end{multline*}
Hence $\phi$ can be factorized through the ideal $J_F$. In such a way we
obtain the desired map $\phi\colon\cal{A}_F
\rightarrow\cal{A}_F\otimes\cal{A}_F$.
\end{pf}

Let $\tau\colon
\cal{A}_F\otimes\cal{A}_F\rightarrow\cal{A}_F\otimes\cal{A}_F$ be a
(involutive) brading given by
$$
\tau(a\otimes b)=(-1)^{\partial a\partial b}\psi(a\otimes b)
$$
where the grading is induced from $W^\wedge$. We have
\begin{gather}
\tau(m_F\otimes\id)=(\id\otimes m_F)(\tau\otimes\id)
(\id\otimes\tau)\label{m-T1}\\
\tau(\id\otimes m_F)=(m_F\otimes\id)(\id\otimes\tau)
(\tau\otimes\id)\label{m-T2}
\end{gather}
as directly follows from \eqref{ps-m1}--\eqref{ps-m2} and
\eqref{ffun}--\eqref{prodcli}.
\begin{lem}
Operators $\tau$ and $\sigma_F$ are mutually compatible,
in the braided sense.
\end{lem}
\begin{pf}
Using \eqref{Sxy} and \eqref{ffun}
we find that the ``braided compatibility''
between $\tau$ and $\sigma_F$ holds on the space $W\otimes W\otimes W$.
Applying inductively \eqref{m-S1}--\eqref{m-S2}
and \eqref{m-T1}--\eqref{m-T2} we conclude that the braided compatibility
is preserved if we pass to the higher filtrant spaces.
\end{pf}

Operators $\tau$ and $\sigma_F$ express twisting properties of the coproduct
map
\begin{lem}
The following identities hold
\begin{gather*}
(\phi\otimes\id)\sigma_F=(\id\otimes\tau)(\sigma_F\otimes\id)(\id\otimes
\phi)=(\id\otimes\sigma_F)(\tau\otimes\id)(\id\otimes\phi)\\
(\phi\otimes\id)\tau=(\id\otimes\tau)(\tau\otimes\id)(\id\otimes\phi)
\phantom{.}\\
(\id\otimes\phi)\sigma_F=(\tau\otimes\id)(\id\otimes\sigma_F)(\phi\otimes\id)
=(\sigma_F\otimes\id)(\id\otimes\tau)(\phi\otimes\id)\\
(\id\otimes\phi)\tau=(\tau\otimes\id)(\id\otimes\tau)(\phi\otimes\id).
\end{gather*}
\end{lem}
\begin{pf}
All these equalities trivially hold on $W\otimes W$. Inductively applying
\eqref{m-T1}--\eqref{m-T2} and \eqref{m-S1}--\eqref{m-S2}, the braided
compatibility between $\tau$ and $\sigma_F$ and the multiplicativity of
$\phi$ we conclude that the above equalities hold on the whole $\cal{A}_F
\otimes\cal{A}_F$.
\end{pf}
As a simple consequence we find
\begin{lem}
The following identity holds
\begin{equation}\label{coas}
(\phi\otimes\id)\phi=(\id\otimes\phi)\phi.
\end{equation}
\end{lem}
\begin{pf}
Let us consider the space $\cal{A}_F\otimes\cal{A}_F\otimes\cal{A}_F$,
endowed with the product $M$ specified by
$$M\leftrightarrow(m_F\otimes m_F\otimes m_F)(\id\otimes
\sigma_F\otimes\sigma_F
\otimes\id)(\id^2\otimes\tau\otimes\id^2).$$
Twisting properties of $\phi$ (and the multiplicativity)
imply that both sides of \eqref{coas} are
unital homomorphisms. On the other hand they trivially coincide on
the generating space $W$. Hence \eqref{coas} holds on the whole $\cal{A}_F$.
\end{pf}
To complete the construction let us introduce the counit and the antipode.
\begin{lem} (i) There exists the unique linear map $\e\colon\cal{A}_F
\rightarrow\Bbb{C}$ satisfying
\begin{align}
\e m_F&=(\e\otimes\e)\sigma_F\tau\label{e-m}\\
\e(W)&=\{0\}\quad\e(1)=1.\label{e-W}
\end{align}

(ii) There exists the unique linear map $\k_F\colon\cal{A}_F\rightarrow
\cal{A}_F$ such that
\begin{gather}
\k_F m_F=m_F(\k_F\otimes\k_F)\tau\sigma_F\tau\sigma_F\tau\label{k-m}\\
\k_F(1)=1\qquad\k_F(x)=-x\label{k-W}
\end{gather}
for each $x\in W$.

(iii) We have
\begin{gather}
m_F(\k_F\otimes\id)\phi=m_F(\id\otimes\k_F)\phi=1\e\label{ant}\\
(\e\otimes\id)\phi=(\id\otimes\e)\phi=\id.\label{cou}
\end{gather}
\end{lem}
\begin{pf}
Formulas of the same form as \eqref{e-m}--\eqref{k-W} consistently and uniquely
determine tensor algebra maps $\e\colon W^\otimes\rightarrow\Bbb{C}$ and
$\k_F\colon W^\otimes\rightarrow W^\otimes$. The consistency follows from
the braided compatibility of $\bigl\{\sigma_F,\tau\bigr\}$ (understood
as braidings on $W^\otimes$). It is therefore sufficient to prove that
$\e(J_F)=\{0\}$ and $\k_F(J_F)=J_F$. In view of the braided-covariance of
the generator relations space $J_F^2$ it is sufficient to check that
$\e(J_F^2)=\{0\}$ and $\k_F(J_F^2)=J_F^2$. A direct computation gives
$$
\e(Q)=0\qquad\k_F(Q)=-Q.
$$
Hence, maps $\e$ and $\k_F$ are factorizable through $J_F$. From the
definition of $\e$ it follows that
\begin{align*}
(\e\otimes\id)\tau&=\id\otimes\e\\
(\id\otimes\e)\tau&=\e\otimes\id.
\end{align*}
Finally, let
us check the antipode and the counit axioms. It
is evident that these equalities hold on the space $W$. On the other
hand, from \eqref{e-m}, \eqref{k-m}, the multiplicativity of $\phi$ and
the twisting properties of $\phi$ and $\e$ it follows that elements
on which equalities \eqref{ant}--\eqref{cou} hold form a
subalgebra of $\cal{A}_F$. Hence \eqref{ant}--\eqref{cou} hold on the whole
$\cal{A}_F$.
\end{pf}

Combining all derived properties we conclude that
$G=\bigl(\cal{A}_F,\bigl\{\phi,\e,\k_F,\sigma_F\bigr\}\bigr)$ is a
braided quantum group.

The presented construction can be performed in a conceptually
different way, explicitly working in terms of the exterior algebra
$W^\wedge$.

Let us consider the space $W\oplus W$, endowed with a braid operator
$\Psi$, given by the block matrix
\begin{equation*}
\Psi=\begin{pmatrix}
\psi&0&0&0\\
0&0&\psi&0\\
0&\psi&0&0\\
0&0&0&\psi
\end{pmatrix}
\end{equation*}
and with the quadratic form $E$ defined as follows. Its restriction on
``homogeneous'' summands in the tensor square of $W\otimes W$ coincides
with $F$, while $E=-F/2$ on ``mixed'' terms.
Let $\Delta\colon W\rightarrow W\oplus W$ be the standard diagonal map.
\begin{lem}\label{lem:7}
There exists the unique algebra map $\Delta^*\colon
\cal{A}_F\rightarrow\cl(W\oplus W,\Psi,E)$ extending $\Delta$. The map
$\Delta^*$ is a homomorphism of the corresponding exterior algebras,
too.
\end{lem}
\begin{pf}
By definition, $\Delta$ intertwines braids $\psi$ and $\Psi$. This
implies that $\Delta$ is uniquely extendible to a unital homomorphism
of the form $\Delta^*\colon W^\wedge\rightarrow(W\oplus W)^\wedge$.
The same map is also a homomorphism of associated Clifford algebras,
as follows from the equality
$$E\bigl(\Delta(x),\Delta(y)\bigr)=F(x,y)$$
and expression \eqref{prodcli} for the Clifford product.
\end{pf}

Algebra $\cal{A}_F$
can be embedded in $\cl(W\otimes W,\Psi,E)$ in two different natural ways,
with the help of monomorphisms $\imb_\pm\colon\cal{A}_F\rightarrow
\cl(W\otimes W,\Psi,E)$ specified by
$$\imb_-(x)=(x,0)\qquad\imb_+(x)=(0,x)$$
where $x\in W$. The maps $\imb_\pm$ are also embeddings of corresponding
exterior algebras. The formula
$$\imb_F(a\otimes b)=\imb_-(a)\imb_+(b)$$
defines a linear map $\imb_F\colon\cal{A}_F\otimes\cal{A}_F\rightarrow
\cl(W\otimes W,\Psi,E)$.
\begin{lem}
The map $\imb_F$ is an algebra isomorphism.
\end{lem}
\begin{pf}It is sufficient to observe that
$$\imb_+(x)\imb_-(y)=-\sum_k\imb_-(y_k)\imb_+(x_k)-F(x,y)1$$
where $x,y\in W$ and $\Sum_ky_k\otimes x_k=\psi(x\otimes y)$.
This gives all mutual
relations between elements
of two embedded subalgebras.
\end{pf}
We have
\begin{equation}\label{d=jf}\Delta^*=\imb_F\phi.\end{equation}
Indeed, both sides of \eqref{d=jf} are unital homomorphisms trivially
coinciding on $W$, by construction.
\section{Concluding Remarks}
Let us analyze a special class of braided Clifford algebras, admitting
spinor representations of a classical type.
Geometrically these algebras represent
particularly regular
quantum spaces, possessing ``sufficiently large'' symmetry groups.
Let us assume that $F$ is nondegenerate, and that
$$W=W_-\oplus W_+$$
where $W_\pm$ are $F$-isotropic subspaces (mutually naturally dual).
Further, let us assume that the above decomposition is compatible with
the braid $\psi$, in the sense that
$$\psi(W_i\otimes W_j)=W_j\otimes W_i$$
for $i,j\in\{+,-\}$.

The corresponding exterior algebras
$W_\pm^\wedge$ are understandable as subalgebras of $\cal{A}_F$, in a
natural manner. Moreover, the map $\mu_F\colon W_-^\wedge\otimes
W_+^\wedge\rightarrow\cal{A}_F$ given by
$$\mu_F(a\otimes b)=ab$$
is bijective. The corresponding ``spinor space'' can be defined as follows
\cite{MZ}.
The restriction $\e_+=\e{\restr}{W_+^\wedge}$
is multiplicative. In particular, $\e_+$ gives a left
$W^\wedge_+$-module structure on $\Bbb{C}$.
On the other hand, $\cal{A}_F$ is a right ${W^\wedge_+}$-module,
in a natural manner. Let ${\cal S}$ be a left $\cal{A}_F$-module, given by
$$ {\cal S}=\cal{A}_F\otimes_*\Bbb{C}$$
where the tensor product is taken over $W_+^\wedge$.

The ``spinor
module'' $\cal{S}$ is simple, and the action of $\cal{A}_F$ is
faithful. In particular, if $\cal{A}_F$ is finite-dimensional (this is a
property of the initial braid $\psi$) then a natural isomorphism
$\cal{A}_F\leftrightarrow \lin(\cal{S})$ holds. The corresponding
quantum space is completely homogeneous and pointless.

Classical Clifford and Weyl algebras are included in the
formalism in a trivial way.

According to Lemma~\ref{lem:7} the map $\Delta^*$ is $F$-independent. It
turns out that $\phi$, and hence $\e$, are $F$-independent, too.

\end{document}